\documentclass[12pt, a4paper]{article}
\usepackage{graphicx}
\usepackage{fullpage} 
\usepackage[english]{babel}
\usepackage{theorem} 
\usepackage{amssymb}
\usepackage{latexsym}

\newtheorem{theo}{Theorem}[section]
\newtheorem{prop}[theo]{Proposition}
\newtheorem{lem}[theo]{Lemma}

\theorembodyfont{\rm} 
\newtheorem{defi}[theo]{Definition}
\newtheorem{rem}[theo]{Remark}

\newenvironment{demo}{\saut\noindent {\it Proof}.}{\hfill$\Box$\par\saut}

\newcommand{\middl}{{\rm mid}}
\newcommand{\slope}{{\rm slope}}
\newcommand{\rien}[1]{}
\newcommand{\saut}{\par\medskip} 
\newcommand{\eps}{\varepsilon}
\newcommand{\vfi}{\varphi}
\newcommand{\disp}{\displaystyle}
\newcommand{\Int}[1]{{\rm Int}\left(#1\right)}

\newcommand{\diam}[1]{\mbox{\rm diam}\left(#1\right)}

\def\IN{\mathbb{N}}
\def\IR{\mathbb{R}}

\title{Transitive sensitive subsystems for interval~maps
\footnotetext{ Studia Math., {\bf 169}, No. 1, 81-104, 2005.}}
\author{Sylvie Ruette\footnote{The author has been partly
supported by a Marie Curie Fellowship of the European Community programme 
Human Potential under contract n$^{\mbox{\tiny o}}$ HPMF-CT-2002-02026.}}
\date{January 18, 2005}

\begin{document}

\maketitle

\begin{abstract}
We state that for continuous interval maps the existence of a non empty
closed invariant subset which is transitive and 
sensitive to initial conditions is implied by positive topological
entropy and implies chaos in the sense of Li-Yorke, and we exhibit
examples showing that these three notions are distinct.
\end{abstract}

\section{Introduction}
In this paper an {\em interval map} is a topological
dynamical system given by a continuous map $f\colon I\to I$ where $I$
is a compact interval. In the literature much has been said about chaos
for interval maps. The point is that the relations between the various
properties related to chaos are much more numerous for these systems than
for general dynamical systems. As a consequence there is a rather ordered
``scale of chaos'' on the interval. For example, for interval maps
topological weak mixing and topological strong mixing are equivalent
\cite{BM}, and transitivity implies sensitivity to initial conditions 
\cite{BM2}, which in turn implies positive topological entropy \cite{Blo3}.
 For more
details on this topic see e.g. \cite{BCop},
\cite[\S\S 6-9]{KS} and \cite{R3}.

Among the different definitions of chaos, a well known one is
{\em chaos in the sense of Li-Yorke}. It follows the ideas of \cite{LY} 
but was formalised later.

\begin{defi}
Let $T\colon X\to X$ be a continuous map on the metric space $X$,
$d$ denoting the distance. The map $T$ is said {\em chaotic
in the sense of Li-Yorke} if there exists an uncountable set $S\subset X$
such that, for all $x,y\in S$, $x\not=y$, one has
$$
\limsup_{n\to+\infty}d(T^n(x),T^n(y))>0\quad\mbox{and}\quad
\liminf_{n\to+\infty}d(T^n(x),T^n(y))=0.
$$
\end{defi}

Note that in the definition of chaos in the sense of Li-Yorke
some people make the extra assumption that
for all $x\in S$ and all periodic points $z\in X$ one has
$\limsup_{n\to+\infty}d(T^n(x),T^n(z))>0$. This gives an equivalent
definition since this property is satisfied by all but at most one points
of the set $S$ \cite[p 144]{BCop}.

Li and Yorke showed that an interval map with a periodic point of
period $3$ is chaotic in the sense of Li-Yorke \cite{LY}. In
\cite{JS} Jankov\'a and Sm\'{\i}tal generalised this result as follows:

\begin{theo}[Jankov\'a-Sm\'{\i}tal]\label{theo:htop-positive-chaos-LY}
If $f\colon I\to I$ is an interval map of positive entropy, then it is
chaotic in the sense of Li-Yorke.
\end{theo}

Recently, Blanchard, Glasner, Kolyada
and Maass proved that, if \mbox{$T\colon X\to X$} is
a continuous map on the compact metric space $X$ such that the topological 
entropy of $T$ is positive, then
the system is chaotic in the sense of Li-Yorke \cite{BGKM}.

The converse of this result is not true, even for interval maps: 
Sm\'{\i}tal \cite{Smi} and Xiong
\cite{Xio3} built interval maps of zero entropy which are chaotic in
the sense of Li-Yorke. See also \cite{MSmi} (a correction is given in 
\cite{Jim}) or \cite{Du2} for examples of a 
$C^{\infty}$ interval map which is  chaotic in the sense of Li-Yorke
and has a null entropy.

\saut
Recall that 
the map $T\colon X\to X$ is {\em transitive} if for all non empty open subsets
$U,V$ there exists an integer $n\geq 0$ such that $T^{-n}(U)\cap V\not=
\emptyset$; if $X$ is compact with no isolated point, $T$ is transitive
if and only if there exists $x\in X$ such that $\omega(x,T)=X$ (where 
$\omega(x,T)$ is the set of limit points of $\{T^n(x)\mid n\geq 0\}$).
The map $T$ has {\em sensitive dependence to initial conditions} 
(or simply
is {\em sensitive}) if there exists $\delta>0$ such that for all $x\in X$ and
all neighbourhoods $U$ of $x$ there exist $y\in U$ and $n\geq 0$ such that 
$d(T^n(x),T^n(y))\geq \delta$. A subset $Y\subset X$ is {\em invariant}
if $T(Y)\subset Y$. 

\medskip
The work of Wiggins \cite{Wig} leads to the following definition
(see, e.g., \cite{FD}).

\begin{defi}
Let $X$ be a metric space. The continuous map $T\colon X\to X$
is said {\em chaotic in the sense of Wiggins} if there exists a non empty
closed invariant subset $Y$ such that the restriction $T|_Y$ is transitive
and sensitive.
\end{defi}

The aim of this paper is to locate this notion with 
respect to the other definitions of chaos.

\begin{rem}
A continuous map $T\colon X\to X$ which is transitive and sensitive is 
sometimes called {\em chaotic in the sense of Auslander-Yorke} \cite{AY}.
If in addition the periodic points are dense, then it is called
{\em chaotic in the sense of Devaney} \cite{Dev}.
\end{rem}

Transitive sensitive subsystems appear naturally when considering a horseshoe,
that is, two disjoint closed intervals $J,K$ such that
$f(J)\cap f(K)\supset J\cup K$, because the points the orbits of which
never escape from $J\cup K$ form a subset on which $f$ acts almost like
a $2$-shift \cite{Bloc3}. For interval maps, positive entropy is equivalent 
to the 
existence of a horseshoe for some power of $f$ \cite{Mis2,BGMY} 
(see also \cite[chap. VIII]{BCop}), thus one can deduce that
a positive entropy interval map has a transitive, sensitive subsystem.
More precisely, Shihai Li proved the following result \cite{Li}.

\begin{theo}[Shihai Li]\label{theo:htop-wiggins}
Let $f\colon I\to I$ be an interval map.
The topological entropy of $f$ is positive if and only if
there exists a non empty closed invariant 
subset $X\subset I$ such that $f|_X$ is transitive, sensitive 
to initial conditions and the periodic points are dense in $X$
(in other words, $f|_X$ is Devaney chaotic).
\end{theo}

In the ``if'' part of this theorem one cannot suppress the assumption on
the periodic points. In Section~\ref{sec:wiggins-htop0} we build a 
counter-example, which leads to the following theorem.

\begin{theo}\label{theo:wiggins-htop0}
There exists a continuous map $f\colon [0,1]\to [0,1]$ 
of zero topological entropy which is chaotic in the sense of Wiggins.
\end{theo}

In \cite{Smi} Sm\'{\i}tal built a zero entropy map $f$ which is chaotic in the
sense of Li-Yorke. If one looks at the construction of $f$, it is not
hard to prove that $f|_{\omega(0,f)}$ is transitive and sensitive to 
initial conditions. We show next theorem in Section~\ref{sec:wiggins-LY}.

\begin{theo}\label{theo:Wiggins-implies-LY}
Let $f\colon I\to I$ be an interval map.
If $f$ is Wiggins chaotic then it is Li-Yorke chaotic.
\end{theo}

The converse of this theorem is not true, contrary to what one may expect
by considering Sm\'{\i}tal's example. The last and longest section is
devoted to the construction of a counter-example that proves the
following result.

\begin{theo}\label{theo:LY-not-wiggins}
There exists a continuous interval map $g\colon I\to I$ which is
chaotic in the sense of Li-Yorke but not in the sense of Wiggins.
\end{theo}

From Theorems \ref{theo:htop-wiggins}, \ref{theo:wiggins-htop0},
\ref{theo:Wiggins-implies-LY}, \ref{theo:LY-not-wiggins} it follows 
that, for interval maps, chaos in the sense of Wiggins  is a strictly
intermediate notion between positive entropy and chaos in the sense of
Li-Yorke. 

\saut
Furthermore the examples of Sections \ref{sec:wiggins-htop0} and
\ref{sec:LY-not-Wiggins} show that the behaviours of zero entropy interval
maps are more varied that one might expect. Let us expose the different
kinds of dynamics exhibited by these maps.

The next result is well known (see, e.g., \cite[p218]{BCop}).

\begin{theo}\label{theo:htop-power-of-2}
Let $f\colon I\to I$ be an interval map. The following properties
are equivalent:
\begin{itemize}
\item the topological entropy of $f$ is zero,
\item every periodic point has a period equal to $2^n$ for some integer
$n\geq 0$.
\end{itemize}
\end{theo}

According to Sharkovskii's Theorem \cite{Sha} the set of periods of periodic
points of a zero entropy interval map is either $\{2^k; 0\leq k\leq n\}$
for some integer $n$ and $f$ is said of type $2^n$, or
$\{2^k;k\geq 0\}$ and $f$ is said of type $2^{\infty}$.
There is little to say about the dynamics of type $2^n$, and some 
interval maps of type $2^{\infty}$ share almost the same dynamics \cite{Del}:
every orbit converges to some periodic orbit of period $2^k$;
these maps are never Li-Yorke chaotic. 

The interval maps of type $2^{\infty}$
that admit an infinite $\omega$-limit set may be Li-Yorke chaotic or
not, as shown by Sm\'{\i}tal \cite{Smi}. A map $f$ that is not Li-Yorke chaotic
is called ``uniformly non-chaotic'' in \cite{BCop} and it satisfies the
following property: every point $x$ is approximately periodic, that is,
for every $\eps>0$ there exists a periodic point $y$ and an integer $N$
such that $|f^n(x)-f^n(y)[<\eps$ for all $n\geq N$.

The maps built in Sections \ref{sec:wiggins-htop0} and
\ref{sec:LY-not-Wiggins} are both zero entropy and Li-Yorke 
chaotic. In the first example there is a transitive sensitive subsystem
which is the core of the dynamics; in particular Li-Yorke chaos can be
read on this subsystem. In the second example this situation does
not occur since there is no transitive sensitive subsystem.

\section{Wiggins chaos implies Li-Yorke chaos}
\label{sec:wiggins-LY}

The following notion of $f$-non separable points was introduced by
Sm\'{\i}tal to give an equivalent condition for chaos in the sense of Li-Yorke
\cite{Smi}. Note that Theorem~\ref{theo:htop0-chaos-LY} was proven to
remain valid for all interval maps by Jankov\'a and Sm\'{\i}tal~\cite{JS}.

\begin{defi}
Let $f\colon I\to I$ be an interval map and $a_0,a_1$ two distinct
points in $I$. The points $a_0, a_1$ are called {\em $f$-separable} if
there exist two disjoint subintervals $J_0,J_1$ and two integers
$n_0,n_1$, such that for $i=0,1$, $a_i\in J_i$, $f^{n_i}(J_i)=J_i$ and
$(f^k(J_i))_{0\leq k<n_i}$ are disjoint. Otherwise they are said {\em
$f$-non separable}.
\end{defi}

\begin{theo}[Sm\'{\i}tal]\label{theo:htop0-chaos-LY}
Let $f\colon I\to I$ be an interval map of zero entropy. The following
properties are equivalent:
\begin{itemize}
\item $f$ is chaotic in the sense of Li-Yorke,
\item there exists $x_0\in I$ such that the set $\omega(x_0,f)$ is
infinite and contains two $f$-non separable points.
\end{itemize}
\end{theo}

In the proof of this theorem, 
Sm\'{\i}tal showed the following intermediate
result which describes the structure of an infinite $\omega$-limit set
of a zero entropy map.

\begin{lem}\label{lem:htop0-Lki}
Let $f\colon I\to I$ be an interval map of zero entropy and $x_0\in I$
such that $\omega(x_0,f)$ is infinite. For all $n\geq 0$
and $0\leq i<2^n$, define 
$$
I_n^i=[\min \omega(f^i(x_0),f^{2^n}), \max\omega(f^i(x_0),f^{2^n})]
\quad\mbox{and}\quad
L_n^i=\bigcup_{k\geq 0}f^{k2^n}(I_n^i).
$$ 
Then
$f(L_k^i)=L_k^{i+1\bmod 2^k}$ for all $0\leq i< 2^k$,
and the intervals $(L_k^i)_{0\leq i< 2^k}$ are pairwise disjoint.
\end{lem}

\begin{lem}\label{lem:fixed-point}
Let $f\colon [a,b]\to \IR$ be a continuous map. If 
$f([a,b])\supset [a,b]$, then $f$ has a fixed point.
\end{lem}

\begin{demo}
There exist $x,y\in [a,b]$ such that
$f(x)\leq a$ and $f(y)\geq b$. One has then
$f(x)-x\leq a-x\leq 0$ and $f(y)-y\geq b-y\geq 0$,
thus there is a point $c\in [x,y]$ such that $f(c)-c=0$.
\end{demo}

\begin{lem}\label{lem:periodic-interval-power-of-2}
Let $f\colon I\to I$ be an interval map of zero entropy. If $J\subset I$ is a
(non necessarily closed) subinterval such that $f^p(J)=J$ and 
$(f^i(J))_{0\leq i<p}$ are pairwise disjoint then $p$ is a power of $2$.
\end{lem}

\begin{demo}
If $J$ is reduced to one point then it is a periodic orbit and by
Theorem~\ref{theo:htop-power-of-2} $p$ is a power of $2$. We assume that
$J$ is non degenerate, which implies that
$f^n(J)$ is a non degenerate interval for all $n\geq 0$.

One has $f^p(\overline{J})=\overline{J}$ thus by Lemma~\ref{lem:fixed-point}
there exists $x\in \overline{J}$ such that $f^p(x)=x$. According to
Theorem~\ref{theo:htop-power-of-2} the period of $x$ is equal to $2^k$
for some $k$; write $p=m2^k$. If $x\in J$ then $(f^i(x))_{0\leq i<p}$ are 
distinct and $p=2^k$. 

Suppose that $m\geq 3$. Then $x\in \partial J$; we assume  
that $x=\sup J$, the case with $x=\inf J$ being
symmetric. One has $x=f^{2^k}(x)\in f^{2^k}(\overline{J})$ and
$f^{2^k}(J)\cap J=\emptyset$ thus $x=\inf f^{2^k}(J)$. But one also has
$x\in f^{2^{k+1}}(\overline{J})$, which contradicts the fact that
$J,f^{2^k}(J), f^{2^{k+1}}(J)$ are pairwise disjoint non degenerate intervals. 
Therefore $m=1$ or $2$ and $p$ is a power of $2$.
\end{demo}

The following result is the key tool in the proof of 
Theorem~\ref{theo:Wiggins-implies-LY}. A rather similar result
can be found in a paper of Fedorenko, Sharkovskii and Sm{\'\i}tal~\cite{FSS}.

\begin{lem}\label{lem:LY-sensitive}
Let $f\colon I\to I$ be an interval map of zero entropy and $x_0$ in $I$
such that $\omega(x_0,f)$ is infinite and does not contain two $f$-non 
separable points. Then
for all $\eps>0$ there exists $\delta>0$ such that, if $x,y\in
\omega(x_0,f)$, $|x-y|<\delta$, then $|f^n(x)-f^n(y)|<\eps$ for all
$n\geq 0$.
\end{lem}

\begin{demo}
Let $X=\omega(x_0,f)$. For all integers $n\geq 0$
and $0\leq i<2^n$ define $a_n^i=\min\omega(f^i(x_0),f^{2^n})$
and $b_n^i=\max\omega(f^i(x_0),f^{2^n})$. Define $I_n^i$ and $L_n^i$
as in Lemma~\ref{lem:htop0-Lki}; one has $I_n^i=[a_n^i,b_n^i]$.  The
points $a_n^i, b_n^i$ belong to $X$ and
\begin{equation}\label{eq:nested-In}
I_{n+1}^i\cup I_{n+1}^{i+2^n}\subset I_n^i\ \mbox{ for all }\ 0\leq i<2^n.
\end{equation}

Suppose that there exists $\eps>0$ such that 
\begin{equation}
\label{eq:In>eps}
\mbox{for all }n\geq 0\mbox{ there is }0\leq i<2^n\mbox{ with }
|I_n^i|\geq\eps.
\end{equation}
Using Equation~(\ref{eq:nested-In}) we can build a sequence $(i_n)_{n\geq
0}$ such that $I_{n+1}^{i_{n+1}}\subset I_n^{i_n}$ and
$|I_n^{i_n}|\geq \eps$ for all $n\geq 0$. Define $I_{\infty}=\bigcap_{n\geq 1}
I_n^{i_n}$.  It is a decreasing intersection of compact 
intervals thus $I_{\infty}$ is
a closed interval and $|I_{\infty}|\geq \eps$. Write $I_{\infty}=[a,b]$; then
$$
a=\lim_{n\to+\infty} a_n^{i_n}\quad\mbox{and}\quad
b=\lim_{n\to+\infty}b_n^{i_n},
$$
thus $a,b\in X$.
One has $a,b\in L_n^{i_n}$ for all $n\geq 0$.
By Lemma~\ref{lem:htop0-Lki}
the intervals $L_{n+2}^{i_{n+2}}$, $f^{2^n}(L_{n+2}^{i_{n+2}})$, 
$f^{2^{n+1}}(L_{n+2}^{i_{n+2}})$ are pairwise disjoint thus 
$\{a,f^{2^n}(a),f^{2^{n+1}}(a)\}$ are distinct.

One has $a\not=b$ and by assumption
$a,b$ are $f$-separable, thus there exist an interval
$J$ and an integer $p\geq 1$ such that $a\in J$, $b\not\in J$,
$f^p(J)=J$ and $(f^i(J))_{0\leq i<p}$ are pairwise disjoint. 
By Lemma~\ref{lem:periodic-interval-power-of-2} 
$p$ is a power of $2$; write $p=2^k$.

Define the interval $K=L_k^{i_k}\cap J$; $K$ contains $a$ and 
$f^{2^k}(K)\subset K$ because $f^{2^k}(L_k^{i_k})=L_k^{i_k}$ by
Lemma~\ref{lem:htop0-Lki}. Thus $K$ contains $\{a,f^{2^k}(a),
f^{2^{k+1}}(a)\}$. These three points belong to $\omega(x_0,f)$ and
are distinct, so one of them belongs to $\Int{K}$ 
and  there exists an integer $n$ such that
$f^n(x_0)\in K$. We have then
$$
X=\omega(x_0,f)\subset \bigcup_{j=0}^{2^k-1} f^j(\overline{K}).
$$

Let $b'\in X$ such that $f^{2^{k+1}}(b')=b$ and
let $0\leq j<2^k$ such that $b'\in f^j(\overline{K})$. The points
$b', f^{2^k}(b')$ and $f^{2^{k+1}}(b')$ belong to $f^j(\overline{K})$ and
they are distinct (same proof as for $a$) 
thus one of them belongs to $f^j(K)$, which implies that
$b\in f^j(K)$. One has $j\not=0$ because $b\not\in J$ and 
$K\subset J$. But on the
other hand $b\in f^j(L_k^{i_k})\cap L_k^{i_k}$ which is empty by 
Lemma~\ref{lem:htop0-Lki}, thus we get a contradiction.
We deduce that Equation~(\ref{eq:In>eps}) is false.

Let $\eps>0$; the negation of Equation~(\ref{eq:In>eps}) implies that
there exists $n\geq 0$ such
that $|I_n^i|<\eps$ for all $0\leq i<2^n$.
Let $\delta>0$ be the minimal distance
between two distinct intervals among $(I_n^i)_{0\leq i<2^n}$. If
$x,y\in X$ with $|x-y|<\delta$ then there exists $0\leq i<2^n$ such
that $x,y\in I_n^i\cap \omega(x_0,f)=\omega(f^i(x_0),f^{2^n})$,
thus for all $k\geq 0$ 
one has  $f^k(x),f^k(y)\in \omega(f^{i+k}(x_0),f^{2^n})\subset
I_n^{i+k\bmod 2^n}$, which implies that $|f^k(x)-f^k(y)|<\eps$. 
\end{demo}

Now we are ready to prove

\saut\noindent{\bf Theorem \ref{theo:Wiggins-implies-LY} }
Let $f\colon I\to I$ be an interval map.
If $f$ is Wiggins chaotic then it is Li-Yorke chaotic.

\saut
\begin{demo}
We show the result by refutation.
Suppose that $f$ is not chaotic in the sense of Li-Yorke. By
Theorem~\ref{theo:htop-positive-chaos-LY} one has $h_{top}(f)=0$.
Consider a closed invariant subset $Y\subset I$ such that $f|_Y$ is
transitive. If $Y$ is finite or has an isolated point then 
$f|_Y$ is not sensitive. If $Y$ is infinite with no isolated point, 
there exists $x_0\in Y$ such that
$\omega(x_0,f)=Y$. By Theorem~\ref{theo:htop0-chaos-LY} $Y$ does not
contain two $f$-non separable points, thus
by  Lemma~\ref{lem:LY-sensitive} $f|_Y$ is not
sensitive. 
\end{demo}


\section{Wiggins chaos does not imply positive entropy}
\label{sec:wiggins-htop0}

We are going to build an interval map of zero entropy which is chaotic in
the sense of Wiggins. It resembles the maps built by Sm\'{\i}tal
(map $f$ in \cite{Smi}) and Delahaye (map $g$ in \cite{Del}), 
however we give full details because this construction will 
be used as a basis for the next example.

\saut\noindent{\bf Notation. } If $I$ is an interval, let
$\middl(I)$ denote the middle of $I$. If $f$ is a linear map, let
$\slope(f)$ denote its constant slope. We write $\uparrow$ (resp.
$\downarrow$) for ``increasing'' (resp. ``decreasing'').

\saut 
Let $(a_n)_{n\geq 0}$ be an increasing sequence of numbers less
that $1$ such that $a_0=0$. Define
$I_0^1=[a_0,1]$ and, for all $n\geq 1$,
$$
I_n^0=[a_{2n-2},a_{2n-1}],\ L_n=[a_{2n-1},a_{2n}],\ I_n^1=[a_{2n},1].
$$
One has $I_n^0\cup L_n\cup I_n^1=I_{n-1}^1$.
We fix $(a_n)_{n\geq 0}$ such that the lengths of the intervals
$I_n^0, I_n^1$ satisfy:
\begin{itemize}
\item if $n$ is odd,  $|I_n^0|=\frac{1}{3^n}|I_{n-1}^1|$
and $|I_n^1|=\left(1-\frac{2}{3^n}\right)|I_{n-1}^1|$,
\item if $n$ is even,  $|I_n^0|=\left(1-\frac{2}{3^n}\right)|I_{n-1}^1|$ and 
$|I_n^1|=\frac{1}{3^n}|I_{n-1}^1|$.
\end{itemize}
This implies that $|L_n|=\frac{1}{3^n}|I_{n-1}^1|$ for all $n\geq 1$.
Note that $|I_n^1|\to 0$, that is, $\lim_{n\to +\infty} a_n=1$; hence
$\bigcup_{n\geq 1}(I_n^0\cup L_n)=[0,1)$.

For all $n\geq 1$, let $\vfi_n\colon I_n^0\to I_n^1$ be the increasing linear
homeomorphism mapping $I_n^0$ onto $I_n^1$; the slope of $\vfi_n$ is
$\slope(\vfi_n)=\frac{|I_n^1|}{|I_n^0|}$. Define the map $f\colon
[0,1]\to [0,1]$ such that $f$ is continuous on $[0,1)$ and
\begin{eqnarray*}
&& f(x)=\vfi_1^{-1}\circ
\vfi_2^{-1}\circ\cdots\circ\vfi_{n-1}^{-1}\circ \vfi_n(x)\mbox{ for all }
x\in I_n^0,\ n\geq 1,\\ 
&& f|_{L_n} \mbox{ is linear for all }n\geq
1,\\ 
&&f(1)=0.
\end{eqnarray*}
Note that $f|_{I_n^0}$ is linear $\uparrow$. 
We will show below that $f$ is continuous at $1$.

\begin{figure}[htb]
\centerline{
\includegraphics[width=12cm]{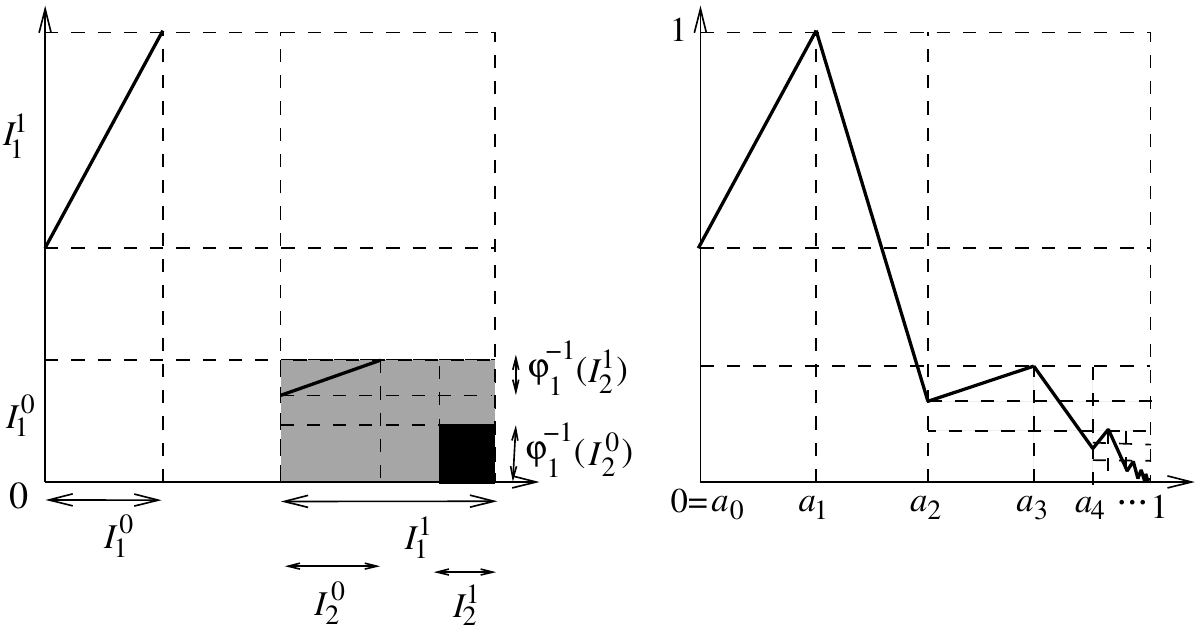}}
\caption{The first steps of the construction of $f$ (left) and
the graph of $f$ (right). This map has a zero entropy and
the invariant set $\omega(0,f)$ is transitive and sensitive. 
\label{fig:wiggins-htop0-partial-construction}}
\end{figure}

Let us explain the underlying construction.  At step $n=1$ the
interval $I_1^0$ is sent linearly onto $I_1^1$  (hence
$f|_{I_1^0}=\vfi_1$) and we decide that $f(I_1^1)\subset I_1^0$ (grey
area on  Figure~\ref{fig:wiggins-htop0-partial-construction}). Then we do the same
kind of construction in the grey area with respect to
$I_2^0,I_2^1\subset I_1^1$: we rescale $I_2^0,I_2^1$ as
$\vfi_1^{-1}(I_2^0),\vfi_1^{-1}(I_2^1) \subset I_1^0$ (on the vertical
axis) then we send linearly $I_2^0$ onto $\vfi_1^{-1}(I_2^1)$; in this
way $f|_{I_2^0}=\vfi_1^{-1}\circ\vfi_2$. We repeat this construction on
$I_2^1$ (black area), and so on. Finally we fill the gaps in a linear
way and we get the whole map, which is
pictured on the right side of
Figure~\ref{fig:wiggins-htop0-partial-construction}. 

\saut
Let $J_0^0=[0,1]$ and for all $n\geq 1$ define the subintervals 
$J_n^0,J_n^1\subset J_{n-1}^0$ by
$\min J_n^0=0$, $\max J_n^1=\max J_{n-1}^0$ and
$\disp \frac{|J_n^i|}{|J_{n-1}^0|}=\frac{|I_n^i|}{|I_{n-1}^1|}$ for $i=0,1$.

\saut
To show that $f$ is continuous at $1$, it is enough to prove that
$\max (f|_{I_n^1})$ tends to $0$ when $n$ goes to infinity. 
For all $n\geq 1$ one has 
\begin{eqnarray*}
\vfi_n(\max I_n^0)&=&\max I_n^1=1=\min I_{n-1}^1+|I_{n-1}^1|\\
\vfi_{n-1}^{-1}\circ \vfi_n(\max I_n^0)&=&
\min I_{n-1}^0+|I_{n-1}^1|\slope (\vfi_{n-1}^{-1})\\
&=&\min I_{n-2}^1+|I_{n-1}^1|\slope (\vfi_{n-1}^{-1})\\
\vfi_{n-2}^{-1}\circ \vfi_{n-1}^{-1}\circ\vfi_n(\max I_n^0)&=& 
\min I_{n-2}^0+|I_{n-1}^1|\slope (\vfi_{n-2}^{-1})\slope (\vfi_{n-1}^{-1})\\
&\vdots &
\end{eqnarray*}

\vspace{-2.8em}
\begin{eqnarray*}
\vfi_1^{-1}\circ \vfi_2^{-1}\circ\cdots\circ 
\vfi_{n-1}^{-1}\circ\vfi_n(\max I_n^0)&=&
\disp \min I_1^0+|I_{n-1}^1|\prod_{i=1}^{n-1}\slope(\vfi_i^{-1})\\
&=&\prod_{i=1}^{n-1}\frac{|I_i^0|}{|I_{i-1}^1|}=|J_{n-1}^0|
\end{eqnarray*}
Consequently,
\begin{equation}\label{eq:wiggins-htop0-max-f(In0)}
f(\max I_n^0)=|J_{n-1}^0|=\max J_{n-1}^0.
\end{equation}
According to the definition of $f$, one has $\max (f|_{I_{n-1}^1})=
f(\max I_n^0)$, thus $\max (f|_{I_{n-1}^1})=|J_{n-1}^0|$.
By definition, $|J_{n-1}^0|\leq\frac{1}{3^{n-2}}$, which tends to $0$,
therefore $f$ is continuous at $1$. 

\saut
Next Lemma describes the action of $f$ on the intervals $(J_n^i)$ and
$(I_n^i)$ and collects the properties that we will use later.
The interval $I_n^1$ is periodic
of period $2^n$ and the map $f^{2^{n-1}}$ swaps $I_n^0$ and $I_n^1$,
however we prefer to deal with $J_n^0=f(I_n^1)$; this will simplify
the proofs because $f|_{I_n^1}$
is not monotone whereas $f^i|_{J_n^0}$ is linear for all 
$1\leq i\leq 2^n-1$.

\begin{lem}\label{lem:wiggins-htop0-summary-Jn0}
Let $f$ be the map defined above. Then, for all
$n\geq 1$,
\begin{enumerate}
\item $f(I_n^1)=J_n^0$,
\item $f(I_n^0)=J_n^1$,
\item $f^i|_{J_n^0}$ is linear $\uparrow$ for all $1\leq i\leq 2^n-1$,
\item $f^{2^{n-1}-1}(J_n^0)=I_n^0$ and $f^{2^n-1}(J_n^0)=I_n^1$,
\item $f^i(J_n^0)\subset \bigcup_{1\leq k\leq n} I_k^0$ for all 
$0\leq i\leq 2^n-2$,
\item $\left(f^i(J_n^0)\right)_{0\leq i<2^n}$ are pairwise disjoint;
\end{enumerate}
and the previous points also imply
\begin{enumerate}
\addtocounter{enumi}{6}
\item $f^{2^{n-1}}(J_n^0)=J_n^1$,
\item $f^{2^n}(J_n^0)=J_n^0$,
\item $f^{2^{n-1}}|_{I_n^0}$ is linear $\uparrow$ and
$f^{2^{n-1}}(I_n^0)=I_n^1$,
\item $f^{2^{n-1}}(I_n^1)=I_n^0$,
\item $(f^i(I_n^0))_{0\leq i<2^n}$ are pairwise disjoint and 
$f^{2^n}(I_n^1)=I_n^1$.
\end{enumerate}
\end{lem}

\begin{demo}
According to Equation~(\ref{eq:wiggins-htop0-max-f(In0)}), 
$\max (f|_{I_n^1})=f(\max I_{n+1}^0)=\max J_n^0$; moreover $f(1)=0=\min
J_n^0$. Thus $f(I_n^1)=J_n^0$ by continuity; this is the point~(i). 

According to the definition of $f$,
\begin{eqnarray*}
|f(I_n^0)|&=&
|I_n^0| \slope(\vfi_n) \prod_{i=1}^{n-1}\slope(\vfi_i^{-1})\\
&= &|I_n^1|\prod_{i=1}^{n-1}\frac{|I_i^0|}{|I_i^1|}
=\frac{|I_n^1|}{|I_{n-1}^1|}\prod_{i=1}^{n-1}\frac{|I_i^0|}{|I_{i-1}^1|}\\
&=&|J_n^1|
\end{eqnarray*}
Moreover $f(\max I_n^0)=\max J_{n-1}^0=\max J_n^1$ according to
Equation~(\ref{eq:wiggins-htop0-max-f(In0)}), 
thus $f(I_n^0)=J_n^1$. This gives the point~(ii).

\saut
We show by induction on $n$ that the points (iii) and (iv) are satisfied.
\begin{itemize}
\item 
This is true for $n=1$ because $J_1^0=I_1^0$, $J_1^1=I_1^1$,
$f|_{I_1^0}=\vfi_1$ is linear $\uparrow$ and $f(I_1^0)=I_1^1$.
\item
Suppose that the induction hypothesis is true for $n$. Since
$J_{n+1}^0\subset J_n^0$, the map $f^i|_{J_{n+1}^0}$ is linear $\uparrow$
for all $1\leq i\leq 2^n-1$ and $f^{2^n-1}(J_{n+1}^0)\subset I_n^1$;
moreover the linearity implies that 
$$
\min f^{2^n-1}(J_{n+1}^0)=\min f^{2^n-1}(J_n^0)=\min I_n^1=\min I_{n+1}^0
$$
and
$$
\frac{|f^{2^n-1}(J_{n+1}^0)|}{|I_n^1|}=\frac{|J_{n+1}^0|}{|J_n^0|}
=\frac{|I_{n+1}^0|}{|I_n^1|}.
$$
Therefore $f^{2^n-1}(J_{n+1}^0)=I_{n+1}^0$.
Then $f^{2^n}(J_{n+1}^0)=J_{n+1}^1$ by the point~(ii). Since $J_{n+1}^1
\subset J_n^0$, the induction hypothesis applies: $f^i|_{J_{n+1}^1}$
is linear $\uparrow$ for all $1\leq i\leq 2^n-1$, $f^{2^n-1}(J_{n+1}^1)
\subset I_n^1$, and by linearity
$$
\max f^{2^n-1}(J_{n+1}^1)=\max f^{2^n-1}(J_n^0)=1=\max I_{n+1}^1
$$
and $|f^{2^n-1}(J_{n+1}^1)|=|I_{n+1}^1|$, hence $f^{2^{n+1}-1}(J_{n+1}^0)\!=\!
f^{2^n-1}(J_{n+1}^1)\!=\!I_{n+1}^1$. This gives the points (iii) and (iv) 
for $n+1$.
\end{itemize}

Now we prove the point (v) by induction on $n$.
\begin{itemize}
\item 
This is true for $n=1$ because $J_1^0=I_1^0$.
\item 
Suppose that the induction hypothesis is true for $n$. Since $J_{n+1}^0
\subset J_n^0$ we have that $f^i(J_{n+1}^0)\subset \bigcup_{1\leq k\leq n}
I_k^n$ for all $0\leq i<2^n-1$. Moreover $f^{2^n-1}(J_{n+1}^0)=
I_{n+1}^0$ by the point~(iv) and $f^{2^n}(J_{n+1}^0)=f(I_{n+1}^0)=
J_{n+1}^1$ by the point~(ii). Since $J_{n+1}^1\subset J_n^0$ we can use
the induction hypothesis again and we get that
$f^{2^n+i}(J_{n+1}^0)\subset \bigcup_{1\leq k\leq n} I_k^n$ for all 
$0\leq i<2^n-1$.
This gives the point~(v) for $n+1$.
\end{itemize}

Next we prove the point~(vi). 
Suppose that $f^i(J_n^0)\cap f^j(J_n^0)\not=\emptyset$ for some
$0\leq i<j<2^n$. Then $f^{2^n-1-j}(f^i(J_n^0))\cap f^{2^n-1-j}(f^j(J_n^0))
\not=\emptyset$. But $f^{2^n-1}(J_n^0)=I_n^1$ by the point~(iv)
and $f^{2^n-1-(j-i)}(J_n^0)\subset [0,\max I_n^0]$ by the point~(v),
thus these two sets are disjoint, which is a contradiction. 
We deduce that $(f^i(J_n^0))_{0\leq i<2^n}$ are pairwise disjoint.

\saut
Finally we indicate how to obtain the other points from the previous ones.
The points (vii) and (viii) are implied respectively by (iv)+(ii) and
(iv)+(i). The point (ix) is implied by (iii)+(iv). The point (x)
is given by (i)+(iv). The point (xi) is given by the combination of
(i), (iv) and (vi).
\end{demo}

Define $K_n=\bigcup_{i\geq 0}f^i(I_n^1)$ for all $n\geq 0$ and $K=\bigcap_{n\geq 0} K_n$.
According to Lemma~\ref{lem:wiggins-htop0-summary-Jn0},
$K_n$ is the disjoint union of the intervals 
$\left(f^i(J_n^0)\right)_{0\leq i\leq 2^n-1}$. The set $K$ has a
Cantor-like construction: at each step a middle part of every connected
component of $K_n$ is removed to get $K_{n+1}$. However $K$ is not a
Cantor set because its interior
is not empty (see Proposition~\ref{prop:wiggins-htop0-omega(0,f)}).
In Proposition~\ref{prop:wiggins-htop0-2-infty} 
we state that the entropy of $f$ is null. Next we show
in Proposition~\ref{prop:wiggins-htop0-omega(0,f)} that the set $\omega(0,f)$
contains $\partial K$. Then we prove that $\omega(0,f)$ is 
transitive and sensitive to initial conditions.

\begin{prop}\label{prop:wiggins-htop0-2-infty}
Let $f$ be the map defined above. Then $h_{top}(f)=0$.
\end{prop}

\begin{demo}
By definition the map $f|_{L_n}$ is linear decreasing thus
$f(L_n)$ is included in $[0,f(\max I_n^0)]$. Moreover $f(\max I_n^0)=\max J_{n-1}^0$
by Equation~(\ref{eq:wiggins-htop0-max-f(In0)}), thus $f(L_n)\subset J_{n-1}^0$.
Then Lemma~\ref{lem:wiggins-htop0-summary-Jn0}(iii) implies  that 
$f^{2^{n-1}}|_{L_n}$ is linear decreasing.
The map $f^{2^{n-1}}|_{I_n^0}$
is linear increasing and $f^{2^{n-1}}(I_n^0)=I_n^1$ by 
Lemma~\ref{lem:wiggins-htop0-summary-Jn0}(ix), 
thus $f^{2^{n-1}}(\min L_n)=\max I_n^1=1$; 
moreover $f^{2^{n-1}}(I_n^1)=I_n^0$ by Lemma~\ref{lem:wiggins-htop0-summary-Jn0}(x),
thus $f^{2^{n-1}}(\max L_n)\in I_n^0$. We deduce that 
$f^{2^{n-1}}(L_n)$ contains $L_n\cup I_n^1$, thus by Lemma~\ref{lem:fixed-point}
there exists $z_n\in L_n$ such that $f^{2^{n-1}}(z_n)=z_n$. The period of
$z_n$ is exactly $2^{n-1}$ because $L_n\subset I_{n-1}^1$ and the
intervals $(f^i(I_{n-1}^1))_{0\leq i<2^n}$ are pairwise disjoint
by Lemma~\ref{lem:wiggins-htop0-summary-Jn0}(xi).

By definition $|L_n|\leq |I_n^1|$ thus
$\slope(f^{2^{n-1}}|_{L_n})\leq -2$. 
If the points $x, f^{2^{n-1}}(x),\ldots, f^{k2^{n-1}}(x)$ belong to $L_n$ then
$|f^{(k+1)2^{n-1}}(x)-z_n|\geq 2^k|x-z_n|$ thus, for all $x\in L_n$,
$x\not=z_n$, there exists $k\geq 1$ such that
$f^{k2^{n-1}}(x)\not\in L_n$. Since $I_{n-1}^1=I_n^0\cup L_n\cup I_n^1$ and
$f^{2^{n-1}}(I_{n-1}^1)=I_{n-1}^1$
by Lemma~\ref{lem:wiggins-htop0-summary-Jn0}(xi),
this implies that $f^{k2^{n-1}}(x)\in I_n^0\cup I_n^1$. In addition
$f^{2^{n-1}}(I_n^0)=I_n^1$ by Lemma~\ref{lem:wiggins-htop0-summary-Jn0}(ix)
thus
$$
\forall x\in I_{n-1}^1,\ x\not=z_n,\ \exists k\geq 0,\
f^k(x)\in I_n^1.
$$
Starting with $I_0^1=[0,1]$, a straightforward induction shows that,
for all $x\in [0,1]$, if the orbit of $x$ does not meet
$\{z_n\mid n\geq 1\}$ then for all integers $n\geq 1$ there exists
$k\geq 0$ such that $f^k(x)\in I_n^1$; in particular 
$\omega(x,f)\subset K$.
According to Lemma~\ref{lem:wiggins-htop0-summary-Jn0}(xi)
the set $K$ contains no periodic point because $K\subset \bigcup_{i\geq 0}f^i(I_n^1)$ for
all $n\geq 1$, thus 
every periodic point is in the orbit of some $z_n$, therefore its
period is a power of $2$.
Finally, $h_{top}(f)=0$ by Theorem~\ref{theo:htop-power-of-2}.
\end{demo}

The orbit of $0$ obviously enters $f^i(J_n^0)$ for all $n\geq 0$ and 
$0\leq i<2^n$, thus $\omega(0,f)$ meets all connected components of $K$.
We show in
next Lemma that $\omega(0,f)$ contains $\partial K$; the proof relies
on the idea that the smaller interval among $J_{n+1}^0$ and $J_{n+1}^1$
contains alternatively either $\min J_n^0$ or $\max J_n^0$ when $n$
varies, so that both endpoints of a connected component of $K$ can be
approximated by small intervals of the form $f^i(J_n^0)$.

\begin{prop}\label{prop:wiggins-htop0-omega(0,f)}
Let $f$ and $K$ be as defined above. Then
$\partial K\subset \omega(0,f)$. In particular $\omega(0,f)$ is
infinite and contains $0$, and $f|_{\omega(0,f)}$ is transitive.
\end{prop}

\begin{demo}
According to the definition of $K$, the connected components of $K$ are
exactly the non empty sets of the form $\bigcap_{n\geq 0}f^{j_n}(J_n^0)$
with $0\leq j_n<2^n$.
Let $y$ be a point in 
$\partial K$. For all $n\geq 0$ there exists $0\leq j_n<2^n$
such that $y\in f^{j_n}(J_n^0)$, and there exists a sequence $(y_n)_{n\geq 0}$
such that $y=\lim_{n\to+\infty}y_n$ and $y_n\in\partial f^{j_n}(J_n^0)=
\{\min f^{j_n}(J_n^0),\max f^{j_n}(J_n^0)\}$. Let $\eps>0$ and $N\geq 0$.
Let $n$ be an even integer such that $\frac{1}{3^{n+1}}<\eps$ and
$|y_n-y|<\eps$, and let $k\geq 0$ such that $k2^{n+1}\geq N$.

Firstly suppose that $y_n=\min f^{j_n}(J_n^0)$. The point $0$ belongs
to $J_{n+1}^0$ and $f^{2^{n+1}}(J_{n+1}^0)=J_{n+1}^0$ by
Lemma~\ref{lem:wiggins-htop0-summary-Jn0}(viii)
thus $f^{k2^{n+1}+j_n}(0)$ belongs to $f^{j_n}(J_{n+1}^0)$. By
Lemma~\ref{lem:wiggins-htop0-summary-Jn0}(iii) one has 
$\min f^{j_n}(J_{n+1}^0)=\min f^{j_n}(J_n^0)=y_n$ and
$$
\frac{|f^{j_n}(J_{n+1}^0)|}{|f^{j_n}(J_n^0)|}=\frac{|J_{n+1}^0|}{|J_n^0|}=
\frac{1}{3^{n+1}}<\eps
$$
thus $|f^{k2^{n+1}+j_n}(0)-y_n|<\eps|f^{j_n}(J_n^0)|\leq \eps$.

Secondly suppose that $y_n=\max f^{j_n}(J_n^0)$. The point
$f^{k2^{n+2}}(0)$ belongs to $J_{n+2}^0$ and $f^{2^{n+1}}(J_{n+2}^0)=
J_{n+2}^1$ by Lemma~\ref{lem:wiggins-htop0-summary-Jn0}(vii) thus
$$
f^{k2^{n+2}+2^{n+1}+2^n+j_n}(0)\in f^{2^n+j_n}(J_{n+2}^1).
$$
According to Lemma~\ref{lem:wiggins-htop0-summary-Jn0}(iii)-(vii) one has
$$
\max f^{2^n+j_n}(J_{n+2}^1)=\max f^{2^n+j_n}(J_{n+1}^0)
=\max f^{j_n}(J_{n+1}^1)
=\max f^{j_n}(J_n^0)=y_n
$$
Moreover
$$
f^{2^n}(J_{n+2}^1)\subset f^{2^n}(J_{n+1}^0)=J_{n+1}^1\subset J_n^0.
$$
Thus
$$
\frac{|f^{j_n+2^n}(J_{n+2}^1)|}{|f^{j_n}(J_n^0)|}=
\frac{|f^{2^n}(J_{n+2}^1)|}{|J_n^0|}
=\frac{|f^{2^n}(J_{n+2}^1)|}{|f^{2^n}(J_{n+1}^0)|}\times
\frac{|J_{n+1}^1|}{|J_n^0|}
=\frac{1}{3^{n+2}}\left(1-\frac{2}{3^{n+1}}\right).
$$
Consequently we get that
$|f^{k2^{n+2}+2^{n+1}+2^n+j_n}(0)-y_n|\leq |f^{j_n+2^n}(J_{n+2}^1)|<\eps$.

In both cases there exists $p\geq N$ such that $|f^p(0)-y_n|<\eps$,
thus $|f^p(0)-y|<2\eps$. This means that $y\in \omega(0,f)$, that is,
$\partial K\subset \omega(0,f)$. The point $\{0\}=\bigcap_{n\geq 0}J_n^0$
belongs to $\partial K$ thus $0\in \omega(0,f)$ and $f|_{\omega(0,f)}$ is
transitive.
Finally, $K_n$ has $2^n$ connected components, each of which containing
$2$ connected components of $K_{n+1}$, thus $K$ has an infinite number
of connected components, which implies that $\partial K$ is infinite.
\end{demo}

In the proof of next proposition, we first show that
$K$ contains a non degenerate connected component $B$.

\begin{prop}\label{prop:wiggins-htop0-non-separable-K}
Let $f$ be the map defined above. Then 
$f|_{\omega(0,f)}$ is sensitive to initial condition.
\end{prop}

\begin{demo}
First we define by induction a sequence of intervals 
$B_n=f^{i_n}(J_n^0)$ for some $0\leq i_n<2^n$ such that $B_n\subset
B_{n-1}$ and $|B_n|=\left(1-\frac{2}{3^n}\right)|B_{n-1}|$ for all $n\geq 1$.
\begin{itemize}
\item 
Take $B_0=J_0=[0,1]$.
\item
Suppose that $B_{n-1}=f^{i_{n-1}}(J_{n-1}^0)$ is already built.
If $n$ is even take $i_n=i_{n-1}$ and $B_n=f^{i_n}(J_n^0)$. The map 
$f^{i_{n-1}}|_{J_{n-1}^0}$ is linear $\uparrow$ by 
Lemma~\ref{lem:wiggins-htop0-summary-Jn0}(iii) and $J_n^0\subset J_{n-1}^0$ thus
$$
\frac{|B_n|}{|B_{n-1}|}=\frac{|J_n^0|}{|J_{n-1}^0|}=1-\frac{2}{3^n}.
$$
If $n$ is odd take $i_n=i_{n-1}+2^{n-1}$ and $B_n=f^{i_n}(J_n^0)$.
According to Lemma~\ref{lem:wiggins-htop0-summary-Jn0}(vii)-(iii) 
one has $B_n=f^{i_{n-1}}(J_n^1)$
and $f^{i_{n-1}}|_{J_{n-1}^0}$ is linear $\uparrow$ thus
$$
\frac{|B_n|}{|B_{n-1}|}=\frac{|J_n^1|}{|J_{n-1}^0|}=1-\frac{2}{3^n}.
$$
\end{itemize}

\saut
Let $B=\bigcap_{n\geq 0}B_n$. This is a compact interval and it is 
non degenerate because
$$
\log |B|=\log |B_0|+\sum_{n\geq 1}\log \left(1-\frac{2}{3^n}\right)
>-\infty.
$$
Moreover $B$ is a connected component of $K$, thus $\partial B\subset
\partial K$. Let $b_0=\min B$ and $b_1=\max B$; one has $b_0, b_1\in
\omega(0,f)$ by Proposition~\ref{prop:wiggins-htop0-omega(0,f)}.

The set $\omega(0,f)$ is included in the periodic orbit of $J_n^0$,
consequently 
$f^k(J_n^0\cap \omega(0,f))=f^k(J_n^0)\cap\omega(0,f)$ for all $k\geq 0$.
Let $\eps>0$ and $k\geq 1$.
There exists $n\geq 0$ such that $|J_n^0|<\eps$ and there exist $x_0,x_1
\in J_n^0\cap \omega(0,f)$ such that $f^{i_n+k2^n}(x_0)=b_0$ and 
$f^{i_n+k2^n}(x_1)=b_1$. Let $\delta=|b_1-b_0|/4$.
The triangular inequality implies that either
$|f^{i_n+k2^n}(0)-b_0)|\geq 2\delta$ or
$|f^{i_n+k2^n}(0)-b_1|\geq 2\delta$. 
In other words, for all $\eps>0$ and $k\geq 1$ there exist $x\in [0,\eps]
\cap \omega(0,f)$ and $i\geq k$ such that $|f^i(0)-f^i(x)|\geq 2\delta$.
Let $y\in \omega(0,f)$ and $\eps>0$; there exists $k\geq 0$ such that
$|f^k(0)-y|<\eps/2$. By continuity of $f^k$ there is $\eta>0$ such that
$f^k([0,\eta])\subset [y-\eps,y+\eps]$. What precedes shows that there
exists $x\in [0,\eta]cap \omega(0,f)$ 
and $i>k$ such that $|f^i(0)-f^i(x)|\geq 2\delta$
thus, if $z=f^k(x)$, $z'=f^k(0)$ and $j=i-n$ we get that $z,z'\in
[y-\eps, y+\eps]$ and $|f^j(z)-f^j(z')|\geq 2\delta$. Then the triangular
inequality implies that either $|f^j(y)-f^j(z)|\geq \delta$ or
$|f^j(y)-f^j('z)|\geq \delta$. We conclude that $f|_{\omega(0,f)}$ is
sensitive to initial conditions.
\end{demo}

Finally Propositions \ref{prop:wiggins-htop0-2-infty}, 
\ref{prop:wiggins-htop0-omega(0,f)} and 
\ref{prop:wiggins-htop0-non-separable-K} give Theorem~\ref{theo:wiggins-htop0}.

\begin{rem}
According to Theorem \ref{theo:Wiggins-implies-LY}  and
\ref{prop:wiggins-htop0-non-separable-K}, the map $f$ is chaotic in the
sense of Li-Yorke. It can be proven directly that $b_0,b_1$ are
$f$-non separable thus Theorem~\ref{theo:htop0-chaos-LY} applies.
\end{rem}

\section{Li-Yorke chaos does not imply Wiggins chaos}\label{sec:LY-not-Wiggins}

The aim of this section is to exhibit an interval map which is
chaotic in the sense of Li-Yorke but has no transitive sensitive subsystem.
This map resembles the one of Section~\ref{sec:wiggins-htop0}:
the construction on the set $\bigcup I_n^0$ is the same except that the
lengths of the intervals differ; the dynamics on $L_n$ is different.

\subsection{Definition of the map {\boldmath $g$}}
\label{subsec:LY-Wiggins-def}
We are going to build a continuous map $g\colon [0,3/2]\to
[0,3/2]$. Let $(a_n)_{n\geq 0}$ be an increasing sequence of numbers
less than $1$ such that $a_0=0$. Define $I_0^1=[a_0,1]$ and for
all $n\geq 1$
$$
I_n^0=[a_{2n-2},a_{2n-1}],\ L_n=[a_{2n-1},a_{2n}],\ I_n^1=[a_{2n},1].
$$
One has $I_n^0\cup L_n\cup I_n^1=I_{n-1}^1$.

\begin{figure}[htb]
\centerline{\includegraphics[width=8cm]{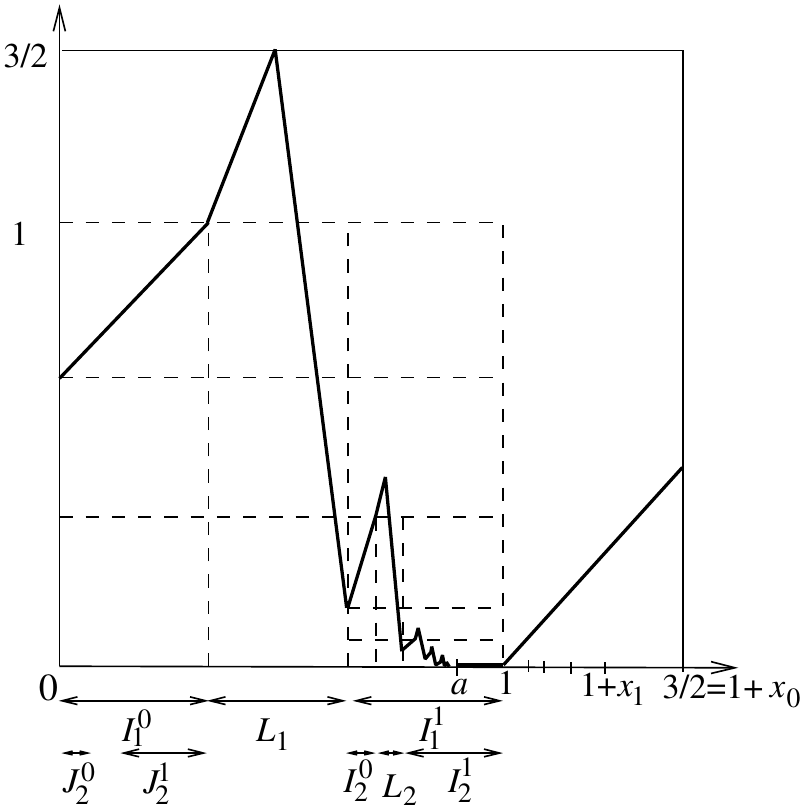}}
\caption{The graph of $g$; this map is Li-Yorke chaotic 
but not Wiggins chaotic. \label{fig:LY-not-wiggins}}
\end{figure}

Fix $(a_n)_{n\geq 0}$ such that the lengths of the intervals
satisfy
$$
\forall n\geq 1,\ |I_n^0|=|L_n|=\frac{1}{3^n}|I_{n-1}^1|\mbox{ and }
|I_n^1|=\left(1-\frac{2}{3^n}\right)|I_{n-1}^1|.
$$
Let $a=\lim_{n\to +\infty}a_n$. 
One has $\bigcup_{n\geq 1}(I_n^0\cup L_n)=[0,a)$ and $a<1$ because
$$
\log (1-a)=\sum_{n=1}^{+\infty}\log \left(1-\frac{2}{3^n}\right)>-\infty.
$$
For all $n\geq 1$,
let $\vfi_n\colon I_n^0\to I_n^1$ be the increasing linear homeomorphism
mapping $I_n^0$ onto $I_n^1$.
Define the map $g\colon [0,3/2]\to  [0,3/2]$ such that $g$ is
continuous on $[0,3/2]\setminus\{a\}$ and
\begin{eqnarray*}
&& g(x)=\vfi_1^{-1}\circ \vfi_2^{-2}\circ\cdots\circ
\vfi_{n-1}^{-1}\circ\vfi_n(x)\mbox{ for all } x\in I_n^0,\ n\geq 1,\\
&& g\mbox{ is linear $\uparrow$ of slope }\lambda_n\mbox{ on } [\min L_n,\middl(L_n)]
\mbox{ for all }n\geq 1,\\
&&g \mbox{ is linear $\downarrow$ on } 
[\middl(L_n),\max L_n]\mbox{ for all }n\geq 1,\\
&&g(x)=0 \mbox{ for all }x\in [a,1],\\
&&g(x)=x-1 \mbox{ for all } x\in [1,3/2],
\end{eqnarray*}
where the slopes $(\lambda_n)$ will be defined below. We will also show below
that $g$ is continuous at $a$. The map $g$ is pictured on 
Figure~\ref{fig:LY-not-wiggins}.

Let $J_0^0=[0,1]$ and for all $n\geq 1$ define the subintervals $J_n^0,
J_n^1\subset J_{n-1}^0$ such that
$\min J_n^0=0$, $\max J_n^1=\max J_{n-1}^0$ and
$\frac{|J_n^i|}{|J_{n-1}^0|}=\frac{|I_n^i|}{|I_{n-1}^1|}\mbox{ for }i=0,1$;
let $M_n=[\max J_n^0,\min J_n^1]$.

Note that on the set $\bigcup_{n\geq 1}I_n^0$ the map $g$ is defined
similarly to the map $f$ in Section~\ref{sec:wiggins-htop0}, thus 
the assertions of Lemma~\ref{lem:wiggins-htop0-summary-Jn0} remain valid for 
$g$, except the point (i) and its derived results (viii), (x), (xi).

\begin{lem}\label{lem:wiggins-LY-summary-Jn0}
Let $g$ be the map defined above.
Then for all $n\geq 1$ one has
\begin{enumerate}
\item $g(I_n^0)=J_n^1$,
\item $g^i|_{J_n^0}$ is linear $\uparrow$ for all $0\leq i\leq 2^n-1$,
\item $g^{2^{n-1}-1}(J_n^0)=I_n^0$ and $g^{2^n-1}(J_n^0)=I_n^1$,
\item $g^i(J_n^0)\subset \bigcup_{1\leq k\leq n}I_k^0$ for all
$0\leq i\leq 2^n-2$,
\item $(g^i(J_n^0))_{0\leq i<2^n}$ are pairwise disjoint.
\item $g^{2^{n-1}}|_{I_n^0}$ is linear $\uparrow$ and 
$g^{2^{n-1}}(I_n^0)=I_n^1$,
\item $g^{2^{n-1}-1}|_{M_n}$ is linear $\uparrow$ and
  $g^{2^{n-1}-1}(M_n)=L_n$,
\item $g(\min L_n)=\min M_{n-1}$,
\item $g^{2^{n-2}}(\min L_n)=\min L_{n-1}$.
\end{enumerate}
\end{lem}

\begin{demo}
For the points (i) to (vi) see the proof of
Lemma~\ref{lem:wiggins-htop0-summary-Jn0}.

According to the point (ii), the map $g^{2^{n-1}-1}|_{M_n}$ is
linear $\uparrow$ because $M_n$ is included in $J_{n-1}^0$. 
Since $M_n=[\max J_n^0,\min J_n^1]$ and
$L_n=[\max I_n^0,\min I_n^1]$, the points (i), (ii) and (iii) imply that
$g^{2^{n-1}-1}(M_n)=L_n$, which is the point~(vii).

The map $g|_{I_n^0}$ is increasing and $\min L_n=\max I_n^0$ thus, 
according to the point~(i), one has $g(\min L_n)=\max J_n^1=\max J_{n-1}^0=
\min M_{n-1}$; this is the point~(viii).

Finally, the points (vii) and (viii) imply the point (ix).
\end{demo}

For all $n\geq 0$, define $x_n=\middl(M_{n+1})$, that is,
$x_n=\frac{3}{2}\prod_{i=1}^{n+1}\frac{1}{3^i}$.
It is a decreasing
sequence and $x_0=1/2$ thus for all $n\geq 0$, $g(1+x_n)$ is well
defined and is equal to $x_n$. 

For all $n\geq 0$ let $t_n=\slope\left(g^{2^n-1}|_{J_n^0}\right)$;
by convention $g^0$ is the identity map so $t_0=1$. Fix $\lambda_1=
\frac{2x_1}{|L_1|}$ and for all $n\geq 2$ define
inductively $\lambda_n$ such that
\begin{equation}
\label{eq:defi-lambda-n}
\frac{|L_n|}{2} \prod_{i=1}^n \lambda_i\prod_{i=0}^{n-2}t_i=x_n.
\end{equation}
By convention an empty product is equal to $1$,
so way Equation~(\ref{eq:defi-lambda-n}) is satisfied for $n=1$.

The slopes $(\lambda_n)_{n\geq 1}$ have been fixed such that
$g^{2^{n-1}}([\min L_n,\middl(L_n)])=[1,1+x_n]$, as 
proven  in next lemma. This means that under the action of $g^{2^{n-1}}$
the image of $L_n$ falls outside of $[0,1]$ but remains close to $1$.
We also list properties of $g$ on
the intervals $L_n$, $I_n^1$ and $[1,1+x_n]$.

\begin{lem}\label{lem:wiggins-LY-summary2}
Let $g$ be the map defined above. Then one has
\begin{enumerate}
\item 
$g^{2^n}|_{[1,1+x_n]}$ is linear $\uparrow$ and
$g^{2^n}([1,1+x_n])=[\min I_{n+1}^0,\middl(L_{n+1})]$ for all $n\geq 0$,
\item 
$g^{2^{n-1}}|_{[\min L_n,\middl(L_n)]}$ is linear $\uparrow$
and $g^{2^{n-1}}([\min L_n,\middl(L_n)])\!=\![1,1+x_n]$ for all $n\geq 1$,
\item
$g^{2^{n+1}}|_{[1,1+x_n]}$ is $\uparrow$ and 
$g^{2^{n+1}}([1,1+x_n])=I_{n+1}^1\cup [1,1+x_{n+1}]$ for
all $n\geq 1$,
\item
$g(I_n^1)\subset [0,\middl(M_n)]$ for all $n\geq 1$,
\item
$g^{2^n}([\min I_n^1,1+x_n])\subset [\min I_n^1,1+x_n]$ and
$g^i([\min I_n^1,1+x_n])\subset [0,1]$ for all $1\leq i\leq 2^n-1$, $n\geq 1$.
\end{enumerate}
\end{lem}

\begin{demo}
The map $g|_{[1,1+x_n]}$ is linear $\uparrow$ and $g([1,1+x_n])=
[0,\middl(M_{n+1})]\subset J_n^0$, thus $g^{2^n}|_{[1,1+x_n]}$ is linear
$\uparrow$ by Lemma~\ref{lem:wiggins-LY-summary-Jn0}(ii). Moreover
$g^{2^n-1}(0)=\min I_{n+1}^0$ and $g^{2^n-1}(\middl(M_{n+1}))=
\middl(L_{n+1})$ by Lemma~\ref{lem:wiggins-LY-summary-Jn0}(iii)+(iv),
 which gives the point~(i).

\saut
Before proving the point~(ii) we show some intermediate results.
Let $n\geq 2$ and $2\leq k\leq n$. One has
\begin{eqnarray*}
\lambda_n\ldots \lambda_k\cdot t_{n-2}\ldots t_{k-2}&=&
\frac{\prod_{i=1}^n\lambda_i\prod_{i=0}^{n-2}t_i}{
\prod_{i=1}^{k-1}\lambda_i\prod_{i=0}^{k-3}t_i}\\
&=&\frac{x_n/|L_n|}{x_{k-1}/|L_{k-1}|}\quad\mbox{by 
Equation~(\ref{eq:defi-lambda-n})}\\
&=&\prod_{i=k+1}^{n+1}\frac{1}{3^i}
\prod_{i=k-1}^{n-1}\frac{1}{1-2/3^i}\times\frac{3^n}{3^{k-1}}\\
&=&\frac{1}{3^{n-k+1}}\prod_{i=k-1}^{n-1}\frac{1}{3^i-2}
\end{eqnarray*}
thus
\begin{equation}\label{eq:wiggins-LY-summary2:1}
\lambda_n\ldots \lambda_k\cdot t_{n-2}\ldots t_{k-2} <1.
\end{equation}

\saut
By definition, $g(\middl(L_n))=g(\min L_n)+\lambda_n \frac{|L_n|}{2}$, and
by Equation~(\ref{eq:defi-lambda-n}), 
\begin{eqnarray*}
\lambda_n \frac{|L_n|}{2}
&=&\frac{x_n}{t_{n-2}\prod_{i=1}^{n-1}\lambda_i\prod_{i=0}^{n-3}t_i}\\
&=& \frac{x_n|L_{n-1}|}{2x_{n-1}t_{n-2}}\\
&=&\frac{1}{3^{n+1}}\frac{|M_{n-1}|}{2}\ \mbox{ because }
t_{n-2}=\frac{|L_{n-1}|}{|M_{n-1}|}\mbox{ by
Lemma~\ref{lem:wiggins-LY-summary-Jn0}(vii)}\\
&<&\frac{|M_{n-1}|}{2}
\end{eqnarray*}
Moreover $g(\min L_n)=\min M_{n-1}$ by 
Lemma~\ref{lem:wiggins-LY-summary-Jn0}(viii) thus
\begin{equation}\label{eq:wiggins-LY-summary2:2}
g([\min L_n,\middl(L_n)])\subset [\min M_{n-1},\middl(M_{n-1})]
\mbox{ for all }n\geq 2.
\end{equation}

\saut
We show by induction on $k=n,\ldots,2$ that\\
-- the map
$g^{2^{n-2}+2^{n-3}+\cdots+2^{k-2}}$
is linear $\uparrow$ of slope $\lambda_n\ldots\lambda_k t_{n-2}\ldots t_{k-2}$
on $[\min L_n,\middl(L_n)]$ and
 maps $\min L_n$ to $\min L_{k-1}$,\\
-- $g^i([\min L_n,\middl(L_n)])\subset [0,1]$ for all
$0\leq i\leq 2^{n-2}+2^{n-3}+\cdots+2^{k-2}$.
\begin{itemize}
\item 
By Equation~(\ref{eq:wiggins-LY-summary2:2})
one has $g([\min L_n,\middl(L_n)])\subset M_{n-1}\subset J_{n-2}^0$, thus
the map $g^{2^{n-2}}|_{[\min L_n,\middl(L_n)]}$ is linear $\uparrow$ of slope
$\lambda_n t_{n-2}$.
According to Lemma~\ref{lem:wiggins-LY-summary-Jn0}(ix) one has
$g^{2^{n-2}}(\min L_n)=\min L_{n-1}$. 
Equation~(\ref{eq:wiggins-LY-summary2:2}) and 
Lemma~\ref{lem:wiggins-LY-summary-Jn0}(iii)+(iv) imply that
$g^i([\min L_n,\middl(L_n)])\subset [0,1]$
for all $1\leq i\leq 2^{n-2}$. This is the induction property at rank
$k=n$.
\item
Suppose that the induction property is true for $k$ with $3\leq k\leq n$.
By Equation~(\ref{eq:wiggins-LY-summary2:1}) one has 
$\lambda_n\ldots\lambda_k\cdot t_{n-2}\ldots t_{k-2}\frac{|L_n|}{2}
\leq \frac{|L_{k-1}|}{2}$ thus 
$g^{2^{n-2}+2^{n-3}+\cdots+2^{k-2}}([\min L_n,\middl(L_n)])\subset
[\min L_{k-1},\middl(L_{k-1})]$. The map $g$ is of slope $\lambda_{k-1}$
on this interval, $g(\min L_{k-1})=\min M_{k-2}$ by 
Lemma~\ref{lem:wiggins-LY-summary-Jn0}(viii) and 
$g([\min L_{k-1},\middl(L_{k-1})])\subset
M_{k-2}$ by Equation~(\ref{eq:wiggins-LY-summary2:2}).
Since $M_{k-2}\subset J_{n-1}^0$, the map
$g^{2^{n-2}+2^{n-3}+\cdots+2^{k-2}+2^{k-3}}$ is linear $\uparrow$ of slope
$\lambda_n\ldots\lambda_{k-1}\cdot t_{n-2}\ldots t_{k-3}$ on
$[\min L_n,\middl(L_n)]$, and it maps $\min L_n$ to $\min L_{k-2}$
by Lemma~\ref{lem:wiggins-LY-summary-Jn0}(ix).
Moreover $g^i([\min L_n,\middl(L_n)])\subset [0,1]$ for all 
$0\leq i\leq 2^{n-2}+2^{n-3}+\cdots+2^{k-2}+2^{k-3}$ by
Lemma~\ref{lem:wiggins-LY-summary-Jn0}(iv) and the induction hypothesis.
This is the property at rank $k-1$.
\end{itemize}

For $k=2$ we finally get that
$g^{2^{n-2}+\cdots +2^0}=g^{2^{n-1}-1}$ is linear $\uparrow$ of slope
$\prod_{i=2}^n\lambda_i\prod_{i=0}^{n-2}t_i$ on $[\min L_n,\middl(L_n)]$,
with $g^{2^{n-1}-1}(\min L_n)=\min L_1$ and $g^{2^{n-1}-1}
([\min L_n,\middl(L_n)])\subset [\min L_1,\middl(L_1)]$.
The map $g$ is of slope $\lambda_1$ on this interval thus,
according to the definition of $\lambda_n$, 
the point~(ii) holds for all $n\geq 2$; it trivially holds for $n=1$ too.
The induction property for $k=2$ also gives that
\begin{equation}\label{eq:wiggins-LY-summary2:5}
g^i([\min L_n,\middl(L_n)])\subset [0,1]\mbox{ for all }
0\leq i\leq2^{n-1}-1,\ n\geq 1.
\end{equation}

\saut
The points (i) and (ii) and 
Lemma~\ref{lem:wiggins-LY-summary-Jn0}(vi) imply the point (iii).

\saut
One has $I_n^1=\bigcup_{k\geq n+1}(I_k^0\cup L_k)\cup [a,1]$.
One can see from the definition of $g$ that 
$$
\max\{g(x)\mid x\in I_k^0\cup
L_k\}=g(\middl(L_k)),
$$ 
thus $g(I_k^0\cup L_k)\subset [0,\middl(M_{k-1})]$ 
by Equation~(\ref{eq:wiggins-LY-summary2:2}). Hence
$$
g(I_n^1)\subset [0,\middl(M_n)]=J_n^0\cup [\min M_n,\middl(M_n)],
$$
which is the point~(iv).

According to Lemma~\ref{lem:wiggins-LY-summary-Jn0}(iii)+(vii), one has that
$g^{2^n-1}(J_n^0)=I_n^1$ and $g^{2^{n-1}-1}([\min M_n,\middl(M_n)])=[\min L_n,
\middl(L_n)]$, and by the point~(ii)
$g^{2^{n-1}}([\min L_n,\middl(L_n)])=[1,1+x_n]$.
Combined with the point~(iv) we get that
\begin{equation}\label{eq:g2n}
g^{2^n}(I_n^1)\subset I_n^1\cup [1,1+x_n].
\end{equation}
Moreover $g^i(J_n^0)\subset [0,1]$ for all $0\leq i\leq 2^n-2$ 
and $g^i([\min M_n,\middl(M_n)])\subset [0,1]$ for all $0\leq i\leq 2^{n-1}-2$
according to Lemma~\ref{lem:wiggins-LY-summary-Jn0}(iv). In addition,
$g^{2^{n-1}+i-1}([\min M_n,\middl(M_n)])=g^i([\min L_n, \middl(L_n)])\subset
[0,1]$ for all \mbox{$0\leq i\leq 2^{n-1}-1$} by 
Equation~(\ref{eq:wiggins-LY-summary2:5}).
Therefore 
\begin{equation}\label{eq:wiggins-LY-summary2:8}
g^i(I_n^1)\subset [0,1]\mbox{ for all }0\leq i <2^n. 
\end{equation}

Finally, $g([1,1+x_n])=[0,\middl(M_{n+1})]\subset J_n^0$ and
the point~(i) implies that $g^{2^n}([1,1+x_n])\subset I_n^1$. Combined
with Equation~(\ref{eq:wiggins-LY-summary2:8}) and (\ref{eq:g2n}) and
Lemma~\ref{lem:wiggins-LY-summary-Jn0}(iv), this gives the point~(v).
\end{demo}

Now we show that
$g$ is continuous at point $a$ as claimed at the beginning of the section.

\begin{lem}
The map $g$ defined above is continuous.
\end{lem}

\begin{demo}
We just have to show the continuity at $a$.
It is clear from the definition that $g$ is continuous at $a^+$.
According to Lemma~\ref{lem:wiggins-LY-summary2}(iv) one has
$g(I_n^1)\subset J_{n-1}^0$. This implies that $g$ is continuous 
at $a^-$ because $\disp\lim_{n\to+\infty} \max J_n^0=0$.
\end{demo}

To end this subsection, let us explain the main underlying ideas of the
construction of $g$ by comparing it with the map $f$ built in 
Section~\ref{sec:wiggins-htop0}. The map $g$ and $f$ are similar on
the set $\bigcup_{n\geq 1} I_n^0$ -- which is the core of the dynamics
of $f$ -- the only difference is the length of the intervals.
For $f$ we showed that $K=\bigcap_{n\geq 0}\bigcup_{i=0}^{2^n-1}
f^i(J_n^0)$ has a non degenerate connected component $B$ and it can be proven
that the endpoints of $B$ are 
$f$-non separable. The same remains true for $g$ with
$B=[a,1]=\bigcap_{n\geq 0}I_n^1$ (the fact that $a,1$ are $g$-non
separable will be proven in Proposition~\ref{prop:wiggins-LY-g-LY-chaotic}).
For $f$ we proved that $\partial K\subset \omega(0,f)$ hence
$\partial B\subset \omega(0,f)$;
for $g$ it is not true that $\{a,1\}\subset \omega(0,g)$ because
the orbit of $0$ stays in $[0,a]$.
The construction of $g$ on the intervals $L_n$ allows to approach $1$ from
outside of $[0,1]$: we will see in 
Proposition~\ref{prop:wiggins-LY-g-LY-chaotic} that $\omega(1+x_0,g)$
contains both $a$ and $1$, which implies chaos in the sense of Li-Yorke.
On the other hand the proof showing that $f|_{\omega(0,f)}$ 
is transitive and sensitive 
fails for $g$ because $\omega(0,g)$ does not contain
$\{a,1\}$ and $\omega(1+x_0,g)$ is not transitive.
We will see in Proposition~\ref{prop:wiggins-LY-g-not-wiggins-chaotic}
that $g$ has no transitive sensitive subsystem at all.

\subsection{{\boldmath $g$} is chaotic in the sense of Li-Yorke}

\begin{prop}\label{prop:wiggins-LY-g-LY-chaotic}
Let $g$ be the map defined in 
Section~\ref{subsec:LY-Wiggins-def}. Then the set $\omega(1+x_0,g)$ is
infinite and contains the points $a,1$, which are $g$-non separable.
Consequently the map $g$ is chaotic in the sense of Li-Yorke.
\end{prop}

\begin{demo}
Lemma~\ref{lem:wiggins-LY-summary2}(iii) implies that
$g^{2^{n+1}}(1+x_n)=1+x_{n+1}$ for all $n\geq 0$.
Since $x_n\to 0$ when $n$ goes to infinity, this implies that
$1\in \omega(1+x_0,g)$. 
Moreover Lemma~\ref{lem:wiggins-LY-summary2}(i) implies that 
$g^{2^n}(1)=\min I_{n+1}^0=a_{2n}$ for all $n\geq 1$, 
hence $a\in \omega(1,g)\subset \omega(1+x_0,g)$. 

Suppose that $A_1,A_2$ are two
periodic intervals such that $a\in A_1$ and $1\in A_2$,
and let $p$ be a common multiple of their periods. One has $g(a)=g(1)=0$
thus $g^p(a)=g^p(1)\in A_1\cap A_2$, and $A_1,A_2$ are not disjoint.
This means that $a,1$ are $g$-non separable. 

It is well known that a finite $\omega$-limit set is cyclic.
Therefore, if $y_0,y_1$ are two distinct points in 
a finite $\omega$-set, the degenerate intervals $\{y_0\}$, $\{y_1\}$
are periodic and $y_0,y_1$ are $g$-separable. This implies that
$\omega(1+x_0,g)$ is infinite. 
We deduce that the map $g$ is chaotic in the sense
of Li-Yorke by Theorem~\ref{theo:htop0-chaos-LY}.
\end{demo}

\subsection{{\boldmath $g$} is not chaotic in the sense of Wiggins}

The main result of this subsection is 
Proposition~\ref{prop:wiggins-LY-g-not-wiggins-chaotic} stating that
$g$ has no transitive sensitive subsystem. 
Next lemma is about the location of transitive subsystems.

\begin{lem}\label{lem:wiggins-LY-transitive-set}
Let $g$ be the map defined in 
Section~\ref{subsec:LY-Wiggins-def} and $Y\subset [0,3/2]$ a closed 
invariant subset with no isolated point such that $g|_Y$ is transitive.
Then
\begin{enumerate}
\item $Y\subset [0,a]$,
\item $\disp Y\subset \bigcup_{i=0}^{2^n-1}g^i(J_n^0)$ for all $n\geq 1$,
\item $g^i(J_n^0\cap Y)=g^i(J_n^0)\cap Y=g^{i\bmod 2^n}(J_n^0)\cap Y$ for all 
$i\geq 0$, $n\geq 0$.
\end{enumerate}
\end{lem}

\begin{demo}
By transitivity there exists $y_0\in Y$ such that $\omega(y_0,g)=Y$; 
in particular the set $Y'=\{g^k(y_0),k\geq 0\}$
is dense in $Y$ and $y\in\omega(y,g)$ for all $y\in Y'$.

Let $n\geq 0$. By Lemma~\ref{lem:wiggins-LY-summary2}(iii),
$g^{2^{n+1}}([1,1+x_n])=I_{n+1}^1\cup [1,1+x_{n+1}]$
thus, according to Lemma~\ref{lem:wiggins-LY-summary2}(v), one obtains 
that for all integers $k\geq 1$,
$g^{k2^{n+1}}([1,1+x_n])\subset I_{n+1}^1\cup [1,1+x_{n+1}]$
and $g^i([1,1+x_n])\subset [0,1]$ for all $i>2^{n+1}$, $i\not\in 2^{n+1}\IN$.
This implies that 
$$
\mbox{for all }i\geq 2^{n+1},\ g^i((1+x_{n+1},1+x_n])\subset 
[0,1+x_{n+1}].
$$ 
Consequently there is no $y\in (1,3/2]=\bigcup_{n\geq 0}
(1+x_{n+1},1+x_n]$ such that $y\in \omega(y,g)$, hence 
$Y'\cap (1,3/2]=\emptyset$ and by density $Y\cap (1,3/2]=\emptyset$.

One has $g^{2^n-1}(0)=a_{2n}$ by 
Lemma~\ref{lem:wiggins-LY-summary-Jn0}(ii)+(iii)
thus the point $0$ is not periodic, hence $\forall k\geq 1$,
$g^k(0)\not\in [a,1]$. 
If $y\in (a,1)$ then $g(y)=0$ and $g^k(y)\not\in [a,1]$ for all
$k\geq 1$, which implies that $y\not\in \omega(y,g)$. Consequently,
$Y\cap (a,1)=\emptyset$. We get that $Y\subset [0,a]\cup\{1\}$, and 
$1\not\in Y$ because $Y$ has no isolated point; this gives the point~(i).

\saut
Let $n\geq 1$. One has $\min L_n=\max I_n^0$ and $\max L_n=\min I_{n+1}^0$
thus $g(\min L_n)=\max J_n^1$ and $g(\max L_n)=\min J_{n+1}^1$ by
Lemma~\ref{lem:wiggins-LY-summary-Jn0}(i). Moreover
$g|_{[\min L_n,\middl(L_n)]}$ is $\uparrow$ and 
$g|_{[\middl(L_n),\max L_n]}$ is linear $\downarrow$ thus there exists
$c_n$ in $[\middl(L_n),\max L_n]$ such that $g(c_n)=g(\min L_n)$.

Since $g([c_n,\max L_n])=[\min J_{n+1}^1,\max J_n^1]$ is included
in $J_{n-1}^0$,
the map $g^{2^{n-1}}|_{[c_n,\max L_n]}$ is linear $\downarrow$ by
Lemma~\ref{lem:wiggins-LY-summary-Jn0}(ii). Moreover $M_n\subset 
g([c_n,\max L_n])$ thus $g^{2^{n-1}}([c_n,\max L_n])$ contains $L_n$ by
Lemma~\ref{lem:wiggins-LY-summary-Jn0}(vii). Consequently 
there exists $z_n\in [c_n,\max L_n]$ such that $g^{2^{n-1}}(z_n)=z_n$
(Lemma~\ref{lem:fixed-point}) and 
$\slope(g^{2^{n-1}}|_{[c_n,\max L_n]})\leq -2$.
Then for every $x\in [c_n,\max L_n]$,
$x\not=z_n$, there exists $k\geq 1$ such that $g^{k2^{n-1}}(x)\not\in
[c_n,\max L_n]$. By Lemma~\ref{lem:wiggins-LY-summary2}(v) one has
$g^{2^{n-1}}(I_{n-1}^1\cup [1,1+x_{n-1}])\subset I_{n-1}^1\cup [1,1+x_{n-1}]$
thus
\begin{eqnarray}  
\lefteqn{\forall x\in [c_n,\max L_n], x\not=z_n,}\nonumber\\
&\exists k\geq 1,\ 
g^{k2^{n-1}}(x)\in I_n^0\cup [\min L_n, c_n]\cup I_n^1\cup [1,1+x_{n-1}].&
\label{eq:out-of-[cn,maxLn]}
\end{eqnarray}

We show by induction on $n$ that 
\begin{equation}\label{eq:Y-In1}
\forall n\geq 0,\ Y'\cap I_n^1\not=\emptyset.
\end{equation}
This is true for $n=0$ because $Y\subset [0,1]=I_0^1$ by the point~(i).
Suppose that there exists $y\in Y'\cap I_{n-1}^1$. Write
$I_{n-1}=I_{n-1}^1=I_n^0\cup L_n\cup I_n^1$; to prove that
$Y'\cap I_n^1\not=\emptyset$ we split into four cases.
\begin{itemize}
\item If $y\in I_n^1$ there is nothing to do.
\item If $y\in I_n^0$ then $g^{2^{n-1}}(y)\in I_n^1$ by
Lemma~\ref{lem:wiggins-LY-summary-Jn0}(vi) and $g^{2^{n-1}}(y)\in Y'$.
\item If $y\in [\min L_n,c_n]$ then $g(y)\in g([\min L_n,\middl(L_n)]$
and $g^{2^{n-1}}(y)\in [1,1+x_n]$ by Lemma~\ref{lem:wiggins-LY-summary2}(ii),
which is impossible because $Y\subset [0,a]$ by the point~(i).
\item If $y\in [c_n,\max L_n]$ then $y\not=z_n$ because $Y$ is infinite.
In addition $g^j(y)\in [0,1]$ for all $j\geq 0$ according to the
point~(i). Then
Equation~(\ref{eq:out-of-[cn,maxLn]}) says that there exists $j\geq 1$
such that $g^j(y')$ belongs to $I_n^0\cup [\min L_n,c_n]\cup I_n^1$ and one of the
first three cases applies.
\end{itemize}

One has $g(I_n^1)\subset J_n^0\cup [\min M_n,\middl(M_n)]$ by
Lemma~\ref{lem:wiggins-LY-summary2}(iv) and
$g^{2^n-1}([\min M_n,\middl(M_n)])=g^{2^{n-1}}([\min L_n,\middl(L_n)])
=[1,1+x_n]$ by Lemmas \ref{lem:wiggins-LY-summary-Jn0}(vii) and
\ref{lem:wiggins-LY-summary2}(ii) respectively. Together with the point~(i)
this implies that
\begin{equation}\label{eq:Y-Jn0}
g(Y\cap I_n^1)\subset J_n^0.
\end{equation}
Equations (\ref{eq:Y-In1}) and (\ref{eq:Y-Jn0}) combined with 
Lemma~\ref{lem:wiggins-LY-summary-Jn0}(i)+(iii) imply that
$$
Y\subset \bigcup_{i=0}^{2^n-1}g^i(J_n^0) \mbox{ for all }n\geq 1,
$$
which is the point~(ii); furthermore
$Y\cap g^i(J_n^0)=Y\cap g^{i\bmod 2^n}(J_n^0)$ for all $i\geq 0$.
Since $g(Y)=Y$ it is clear that $g^i(J_n^0\cap Y)\subset 
g^i(J_n^0)\cap Y$ and that $g^{2^n}(g^i(J_n^0)\cap Y)\subset 
g^{2^n+i}(J_n^0)\cap Y$, thus 
$$
g^i(J_n^0\cap Y)=g^i(J_n^0)\cap Y=g^{i\bmod 2^n}(J_n^0)\cap Y \mbox{ for all }
i\geq 0,
$$
which concludes the proof of the Lemma.
\end{demo}

Next lemma is the key tool in the proof of 
Proposition~\ref{prop:wiggins-LY-g-not-wiggins-chaotic}. 
It relies on the knowledge of the precise location of $g^i(J_n^0)$ in
$\bigcup_{1\leq k\leq n}I_n^0$.

\begin{lem}\label{lem:slope}
Let $g$ be the map defined in 
Section~\ref{subsec:LY-Wiggins-def}.
For all $n\geq 1$ and all $0\leq k\leq 2^n-1$ one has
$\slope\!\left(g^{2^n-1-k}|_{g^k(J_n^0)}\right)\geq 1$.
\end{lem}

\begin{demo}
A {\em (finite) word} $B$ is an element of $\IN^n$ for some $n\in \IN$.
If $B,B'$ are two words, $BB'$ denotes their concatenation and
$|B|=n$ is the length of $B$.

We define inductively a sequence of words $(B_n)_{n\geq 1}$ by:
\begin{itemize}
\item $B_1=1$,
\item $B_n=nB_1B_2\ldots B_{n-1}$.
\end{itemize}
and we define the infinite word $\omega=(\omega(i))_{i\geq 1}$ by
concatenating the $B_n$'s: $\omega=B_1B_2B_3\ldots B_n\ldots$
A straightforward induction shows that
$|B_n|=2^{n-1}$ thus $|B_1|+|B_2|+\cdots+|B_k|=2^k-1$ and
the word $B_{k+1}$ begins at $\omega(2^k)$, which gives
\begin{equation}
\label{eq:omega(2k)}
\omega(2^k)=k+1,
\end{equation}
and
\begin{equation}
\label{eq:omega-repeat}
\omega(2^k+1)\ldots\omega(2^{k+1}-1)=B_1\ldots B_k=
\omega(1)\ldots\omega(2^k-1).
\end{equation}

We prove by induction on $k\geq 1$ that
\begin{equation}
\label{eq:gi(Jn)-omega}
g^{i-1}(J_n^0)\subset I_{\omega(i)}^0\mbox{ for all } n\geq k,\ 
1\leq i\leq 2^k-1.
\end{equation}
\begin{itemize}
\item 

Case $k=1$: $J_n^0\subset I_1^0=I_{\omega(1)}^0$ for all $n\geq 1$.
\item
Suppose that Equation~(\ref{eq:gi(Jn)-omega}) holds for $k$ and let
$n\geq k+1$. Since  $J_n^0$ is included in $J_{k+1}^0$, 
Lemma~\ref{lem:wiggins-LY-summary-Jn0}(iii) implies that
$g^{2^k-1}(J_n^0)\subset I_{k+1}^0$, thus
$g^{2^k}(J_n^0)\subset J_k^0$ by 
Lemma~\ref{lem:wiggins-LY-summary-Jn0}(i). By
induction one has $g^{i-1}(J_k^0)\subset I_{\omega(i)}^0$ for all
$1\leq i\leq 1^k-1$, and by Equation~(\ref{eq:omega-repeat}) , one has
$\omega(i)=\omega(2^k+i)$ for all $1\leq i\leq 2^k-1$. Consequently
$g^{2^k+i-1}(J_n^0)\subset I_{\omega(2^k+i)}^0$ for all $1\leq i\leq 2^k-1$.
Together with the induction hypothesis this gives 
Equation~(\ref{eq:gi(Jn)-omega}) for $k+1$.
\end{itemize}

\saut
Let $\mu_n=\slope(g|_{I_n^0})$. By definition of $g$ one has
$$
\mu_n=\frac{\slope(\vfi_n)}{\prod_{i=1}^{n-1}\slope(\vfi_i)}.
$$
It is straightforward from Equation~(\ref{eq:gi(Jn)-omega}) that
for all $2\leq k\leq 2^n-1$
\begin{equation}\label{eq:slope-gk}
\slope (g^{k-1}|_{J_n^0})=\prod_{i=1}^{k-1} \mu_{\omega(i)}.
\end{equation}
By Lemma~\ref{lem:wiggins-LY-summary-Jn0}(ii)+(iii) 
the map $g^{2^n-1}|_{J_n^0}$
is linear and $g^{2^n-1}(J_n^0)=I_n^1$, thus
$$
\slope(g^{2^n-1}|_{J_n^0})=\frac{|I_n^1|}{|J_n^0|}
=\prod_{i=1}^n\frac{1-2/3^i}{1/3^i}.
$$
Since $\slope(\vfi_i)=\frac{|I_i^1|}{|I_i^0|}=\frac{1-2/3^i}{1/3^i}$, we
get
\begin{equation}
\label{eq:slope-g2n-1}
\slope(g^{2^n-1}|_{J_n^0})=\prod_{i=1}^{2^n-1}\mu_{\omega(i)}=
\prod_{i=1}^n \slope(\vfi_i).
\end{equation}

We show by induction on $n\geq 1$ that for all $1\leq k\leq 2^n-1$
\begin{equation}
\label{eq:prod-mui}
\prod_{i=1}^{k}\mu_{\omega(i)}=\prod_{i=1}^n\slope(\vfi_i)^{\eps_i}
\mbox{ for some }\eps_i=\eps(i,k,n)\in\{0,1\}.
\end{equation}
\begin{itemize}
\item
$\mu_{\omega(1)}=\mu_1=\slope(\vfi_1)$; this gives the case $n=1$.
\item
Suppose that the induction hypothesis is true for
$n$. One has
\begin{eqnarray*}
\prod_{i=1}^{2^n}\mu_{\omega(i)}&=&\prod_{i=1}^{2^n-1}\mu_{\omega(i)}\times
\mu_{n+1}\mbox{ by Equation~(\ref{eq:omega(2k)})}\\
&=&\disp \prod_{i=1}^n\slope(\vfi_i)
\frac{\slope(\vfi_{n+1})}{\prod_{i=1}^n\slope(\vfi_i)} 
\mbox{ by Equation~(\ref{eq:slope-g2n-1})}\\
&=&\slope(\vfi_{n+1})
\end{eqnarray*}
This is Equation~(\ref{eq:prod-mui}) for $n+1$ and $k=2^n$ with 
$\eps(i,k,n1)=0$ for $1\leq i\leq n$ and $\eps(n+1,2^n,n+1)=1$.

Next, $\omega(2^n+1)\ldots\omega(2^{n+1}-1)=\omega(1)\ldots
\omega(2^n-1)$ by Equation~(\ref{eq:omega-repeat}) thus, if
$2^n+1\leq k\leq 2^{n+1}-1$ one has
\begin{eqnarray*}
\prod_{i=1}^k\mu_{\omega(i)}&=&
\prod_{i=1}^{2^n}\mu_{\omega(i)}\prod_{i=2^n+1}^k\mu_{\omega(i)}\\
&=&\slope(\vfi_{n+1})\prod_{i=1}^{k-2^n}\mu_{\omega(i)}\\
&=&\slope(\vfi_{n+1})\prod_{i=1}^n\slope(\vfi_i)^{\eps_=(i,k_2^n,n)}
\end{eqnarray*}
That is, Equation~(\ref{eq:prod-mui}) holds with 
$\eps(i,k,n+1)=\eps(i,k-2^n,n)$ for all $1\leq i\leq n$ and 
$\eps(n+1,k,n+1)=1$. This concludes the induction.
\end{itemize}

Equations (\ref{eq:slope-gk}) and (\ref{eq:prod-mui}) show that for all
$1\leq k\leq 2^n-1$
\begin{equation}
\label{eq:slope-gk-bis}
\slope(g^k|_{J_n})=\prod_{i=1}^{k+1}\mu_{\omega(i)}=
\prod_{i=1}^{n}\slope(\vfi_i)^{\eps_i}\mbox{ for some }
\eps_i\in\{0,1\}.
\end{equation}
Since $\slope\!\left(g^{2^n-1-k}|_{g^k(J_n^0)}\right)=
\frac{\slope(g^{2^n-1}|_{J_n^0})}{\slope(g^k|_{J_n^0})}$, Equations 
(\ref{eq:slope-g2n-1}) and (\ref{eq:slope-gk-bis})
imply that $\slope\!\left(g^{2^n-1-k}|_{g^k(J_n)}\right)$ is a
product of at most $n$ terms of the form $\slope(\vfi_i)$. This 
concludes the proof of the
lemma because $\slope(\vfi_i)\geq 1$ for all $i\geq 1$.
\end{demo}

\begin{prop}\label{prop:wiggins-LY-g-not-wiggins-chaotic}
The map $g$ defined in 
Section~\ref{subsec:LY-Wiggins-def} is not Wiggins chaotic.
\end{prop}

\begin{demo}
Let $Y\subset [0,3/2]$ be a closed invariant subset such that $g|_Y$ is
transitive. We assume that $Y$ has no isolated point, otherwise $g|_Y$
is not sensitive. 

The sets $\left(g^i(J_n^0\cap Y)\right)_{0\leq i\leq 2^n-1}$ are closed and 
by Lemma~\ref{lem:wiggins-LY-summary-Jn0}(v) 
they are pairwise disjoint; let $\delta_n>0$
be the minimal distance between two of these sets.
If $x,x'\in Y$ and $|x-x'|<\delta_n$ then there is $0\leq i\leq 2^n-1$ such
that $x,x'\in g^i(J_n^0)$ and for all $k\geq 0$ one has
$g^k(x),g^k(x')\in g^{i+k\bmod{2^n}}(J_n^0)$ by
Lemma~\ref{lem:wiggins-LY-transitive-set}(ii)+(iii). Let 
$$
\eps_n=\max\{\diam{g^i(J_n^0)\cap Y}\mid 0\leq i<2^n\}.
$$
By Lemma~\ref{lem:slope}, we get that $\diam{g^k(J_n^0)\cap Y}
\leq \diam{g^{2^n-1}(J_n^0)\cap Y}$ for all $0\leq k\leq 2^n-1$.
By Lemma~\ref{lem:wiggins-LY-summary-Jn0}(iii) one has 
$ g^{2^n-1}(J_n^0)=I_n^1$ and by
Lemma~\ref{lem:wiggins-LY-transitive-set}(i) one has 
$I_n^1\cap Y\subset [a_{2n},a]$, thus $\eps_n\leq \diam{I_n^1\cap Y}\leq
a-a_{2n}$, which implies that 
$$
\lim_{n\to+\infty}\eps_n=0.
$$ 
This implies that $g|_Y$ is not sensitive. 
\end{demo}

At last this example is completed. 
Theorem~\ref{theo:LY-not-wiggins} is given by
Propositions \ref{prop:wiggins-LY-g-LY-chaotic} and
\ref{prop:wiggins-LY-g-not-wiggins-chaotic}.



\noindent
Laboratoire de Math\'ematiques --
Topologie et Dynamique -- B\^atiment 425~-- Universit\'e Paris-Sud --
F-91405 Orsay cedex -- France.\\
E-mail: {\tt Sylvie.Ruette@math.u-psud.fr}

\end{document}